\documentclass[12pt,leqno]{article}
\usepackage{amsmath,amssymb}

\title{\it TOPOLOGICAL GROUPS WITH SEVERAL DISCONNECTEDNESS}
\author{MASASI HIGASIKAWA
\thanks{Work partly supported by JSPS Research Fellowships for
Young Scientists and a Grant-in-Aid for Scientific Research
from the Ministry of Education, Science and Culture of Japan.}
}
\date{December 2, 2000}

\newenvironment{Enumerate}{\begin{enumerate}

}{\end{enumerate}}

\newtheorem{theorem}{\sc Theorem}[section]
\newtheorem{prop}[theorem]{\sc Proposition}
\newtheorem{remark}[theorem]{\sc Remark}
\newtheorem{prob}{\sc Problem}
\newtheorem{example}[theorem]{\sc Example}
\newtheorem{defn}[theorem]{\sc Definition}

\newenvironment{proof}{Proof.}{$\blacksquare$\\}

\newcommand{\Pair}[1]{\left\langle #1\right\rangle}
\newcommand{\Z}{{\mathbf Z}}
\newcommand{\Q}{{\mathbf Q}}
\newcommand{\R}{{\mathbf R}}
\newcommand{\Sym}{{\rm Sym}}

\begin{document}
\maketitle

\begin{abstract}
We investigate some properties of topological groups related to 
disconnectedness or Archimedeanness. We prove or disprove the preservation of
those under operations as subgroups, quotients, products, etc.
Characterizations of non-Archimedeanness are obtained by using embedding into
universal groups. We also clarify the differences between the properties by
constructing miscellaneous Polish groups.
\end{abstract}

{
\renewcommand{\thefootnote}{}
\footnotetext{2000 {\it Mathematics Subject Classification}: 22A05.}
}

\section{Introduction}

For topological groups, there are several notions of connectedness or
disconnectedness; some are topological and others are topologico-algebraic.
The topological groups treated here may be connected or zero-dimensional at the
extremes as to the former, and Archimedean or non-Archimedean (in the
terminology of \cite{Sh2}, \cite{Sh1}) as to the latter.

Our main concern is for topologico-algebraic properties, i.e., richness or
poorness in open subgroups. We investigate the interrelation among each other
and with some operations and topological properties.

In Section 2, we introduce the properties and obtain preservation results under
subgroups, quotients, direct products, topological closures and group
extensions. Section 3 is for the characterizations of non-Archimedeanness.
We show universal groups in general and Abelian case, respectively. In
Section 4, we mainly treat Abelian Polish groups. We exhibit peculiar examples
with or without several properties.

All topological spaces we consider here are assumed to be Hausdorff unless
otherwise stated. In particular, a quotient of a topological group is meant to
be that by a closed subgroup.

We mainly follow the notation of standard texts such as \cite{HR}. The neutral
element of a group may be denoted by $1$ (or $0$ in Abelian case). We also
adopt some set-theoretic notations. Let ${}^XY$ denote the set of functions
from $X$ into $Y$ possibly with some induced structure. The set $\{0,1,2,...\}$
of natural numbers is identified with the ordinal number $\omega$ and the
cardinal number $\aleph_0$.

\section{Many or few open subgroups}

We consider following properties for topological groups concerning the
magnitude of open subgroups.

\begin{defn}\rm
\begin{Enumerate}
\item A topological group is {\em topologically non-Archimedean} ({\em TNA})
 if every neighborhood of the neutral element contains an open subgroup.
\item A topological group has {\em sufficiently many open subgroups}
({\em SMOG}) if the intersection of all open subgroups consists of only the
neutral element.
\item A topological group is {\em topologically Archimedean} ({\em TA}) if it
has no proper open subgroups.
\end{Enumerate}
\end{defn}

These have close relation to connectedness as open subgroups are clopen. We
recall some previously known or easily observable facts.

\begin{remark}\label{dc}\rm
\begin{Enumerate}
\item The terminology of (non-)Archimedeanness for Abelian topological groups
 appear in \cite{Sh2} and \cite {Sh1}. As noted in \cite{Sh1}, a topological
 group is TNA if and only if it is Archimedean with respect to the right (or
 left) group uniformity. Being TA is equivalent to each of the following
 properties: every non-empty open set (neighborhood of $1$, respectively)
 generates the whole group.
\item The following implications and incompatibility are straightforward:
\begin{center}
\begin{tabular}{ccc}
TNA          & $\Rightarrow$ & zero-dimensional \\
$\Downarrow$ &               & $\Downarrow$     \\
SMOG         & $\Rightarrow$ & totally-disconnected
\end{tabular},
\end{center}
\[\rm connected \Rightarrow TA;\]
A TA group with SMOG is trivial.
\item If a topological group has SMOG, then the open subgroups induce another
 group topology, which is TNA. So a topological group has SMOG if and only if
 it has a weaker (Hausdorff) TNA group topology.
\item For a locally compact group, each arrow above can be reversed
(cf.~\cite{HR}, Theorem 7.7, Corollary 7.9).
\item As remarked in \cite[p.~9]{BK}, the Polish group $\overline{F}$ in
 \cite[pp.~77--78]{Do} is zero-dimensional but not TNA; it has SMOG. Ancel,
 Dobrowolski and Grabowski \cite{ADG} investigate topological groups with
 similar properties in Banach spaces. Our $\Gamma_0$ in Example \ref{ss} is
 typical one of them.
\item A totally-disconnected TA Polish group is constructed by T.C.~Stevens
 \cite{St} (see also Hartman, Mycielski, Rolewicz and Schinzel \cite{HMRS} and
 Schinzel \cite{Sc}), which is not zero-dimensional (Theorem \ref{non0}).
\end{Enumerate}
\end{remark}

We have preservation under some operations. Some are immediate by definition.

\begin{prop}
\begin{Enumerate}
\item Subgroups, quotient groups and products of TNA groups are also TNA.
\item Having SMOG is hereditary to subgroups and products (but not to 
 quotients, see Example \ref{ss}).
\item Quotient groups and products of TA groups are also TA; subgroups need
 not.
\end{Enumerate}
\end{prop}
\begin{proof}
We show (3) for a product. Suppose that $\Pair{G_i:i\in I}$ is a family of TA
groups. Let $J$ be a finite subset of $I$ and $U_i$ a neighborhood of $1$ in
$G_i$ for $i\in J$. We shall prove that the neighborhood
$U=\prod_{i\in I\setminus J}G_i \times \prod_{i\in J}U_i$ generates the whole
product $\prod_{i\in I}G_i$. For any $g=\Pair{g_i:i\in I}\in\prod_{i\in I}G_i$,
set $g'=\Pair{g_i:i\in I\setminus J}\cup\Pair{1_{G_i}:i\in J}$, i.e.,
$g'=\Pair{g'_i:i\in I}$ with $g'_i=g_i$ for $i\in I\setminus J$ and
$g'_i=1_{G_i}$ for $i\in J$. For each $i\in J$, there exists a finite sequence
$g_{i1},...,g_{in_i}$ in $U_i\cup U_i^{-1}$ such that
$g_i=g_{i1}\cdots g_{in_i}$.
We may assume that $U_i^{-1}=U_i$ and $n_i=n$ for all $i\in J$. Let
$h_k$ for $1\leq k\leq n$ denote
$\Pair{1_{G_i}:i\in I\setminus J}\cup\Pair{g_{ik}:i\in J}$. Then
$g',h_1,...,h_n\in U$ and $g'h_1\cdots h_n=g$. So we are done.
\end{proof}

Topological closure preserves two properties of the three.

\begin{prop}
\begin{Enumerate}
\item A topological group is TA if it has a dense TA subgroup.
\item TNA-ness of a dense subgroup implies that of the whole group.
\end{Enumerate}
\end{prop}
\begin{proof}
(2) Assuming that $H$ is a dense TNA subgroup of $G$, we shall show that $G$ is
also TNA. For a neighborhood $U$ of $1$ in $G$, there exists another $V$ such
that $\overline{V}\subseteq U$. By the assumption, we find an open subgroup
$K$ of $H$ with $K\subseteq V\cap H$. Then the closure $\overline{K}$ in $G$ is
an open subgroup contained in $U$.
\end{proof}

\begin{example}\rm
The additive group of $\R\times\Q_p$ has a dense subgroup $\{(x,x):x\in\Q\}$.
While the latter has SMOG, the former does not.
\end{example}

Next we consider group extensions. Suppose that $G$ is a topological group and
$N$ is a normal subgroup.

\begin{prop}
If $N$ and $G/N$ are TA, so is $G$.
\end{prop}
\begin{proof}
Let $H$ be an open subgroup of $G$. Then $H\cap N$ is open in $N$, and hence
$H\cap N=N$, i.e., $N\subseteq H$. Since $H/N$ is open in $G/N$, we also have
$H/N=G/N$. Therefore $H=G$.
\end{proof}

\begin{theorem}
If both $N$ and $G/N$ are TNA, then so is $G$.
\end{theorem}
\begin{proof}
Let $U$ be a neighborhood of $1$ in $G$. We shall find an open subgroup $H$
contained in $U$.

We choose neighborhoods $U_0$, $V$ and $W$ of $1$ in $G$ as follows. First
let $U_0$ be such that $U_0^2\subseteq U$. By the assumption, there is an open
subgroup $M$ of $N$ contained in $N\cap U_0$. Let $V\subseteq U_0$ be open with
$V^{-1}=V$ and $V^3\cap N\subseteq M$. We denote by $\pi$ the natural
homomorphism $G\rightarrow G/N$. Since $\pi(V)$ is open in $G/N$, it contains
an open subgroup $K$. We set $W=V\cap\pi^{-1}(K)$.

We show that $W^2\subseteq WM$. Suppose that $w_0,w_1\in W$. Since
$\pi(w_0),\pi(w_1)\in K$, we have $\pi(w_0w_1)\in K$. So there is $w_2\in W$
with $\pi(w_2)=\pi(w_0w_1)$. Then $w_2^{-1}w_0w_1\in N\cap W^3\subseteq M$,
and hence $w_0w_1\in w_2 M$.

Now let $H$ be the subgroup of $G$ generated by $W$. Then $H$ is open and
$H\subseteq WM \subseteq U_0^2 \subseteq U$ as desired.
\end{proof}

The ``three-group property'' above does not hold for SMOG
(Example \ref{3smog}).

\section{Universal non-Archimedean groups}

TNA groups are characterized by embeddings into certain topological groups.

For a (usually infinite) set $\Omega$, $\Sym(\Omega)$ denotes the symmetric
group on $\Omega$. We topologize $\Sym(\Omega)$ as a subspace of the product
${}^\Omega\Omega$ with $\Omega$ discrete, which makes $\Sym(\Omega)$ a
topological group.

We recall some cardinal invariants for topological spaces (cf.~\cite{Hode}).
Suppose $X$ is a topological space. The {\em weight}, the {\em cellularity} and
the {\em character} of $X$ are defined respectively as follows:
\[ w(X) = \min\{|{\cal B}|:\mbox{$\cal B$ a base for $X$}\}+\aleph_0,\]
\[ c(X) = \sup\{|{\cal V}|:
 \mbox{$\cal V$ a disjoint family of open sets in $X$}\}+\aleph_0,\]
\[ \chi(X) = \sup\{\chi(p,X):p\in X\}+\aleph_0,\]
where
\[ \chi(p,X) = \min\{|{\cal V}|:\mbox{$\cal V$ a local base for $p$}\}.\]
Note that $c(X),\chi(X)\leq w(X)$ and that in case $X$ is a topological group,
\[ \chi(X) = \chi(1,X)+\aleph_0.\]

\begin{theorem}\label{tna}
\begin{Enumerate}
\item A topological group is TNA if and only if it is isomorphic to a subgroup
 of a symmetric group. Specifically, for a topological group $G$ and an
 infinite cardinal $\kappa$, the following are equivalent.
\begin{Enumerate}
\item $G$ is isomorphic to a subgroup of $\Sym(\kappa)$.
\item $G$ is TNA and $w(G)\leq\kappa$.
\item $G$ is TNA and $c(G),\chi(G)\leq\kappa$.
\end{Enumerate}
\item An Abelian topological group is TNA if and only if it is isomorphic to a
 subgroup of a direct product of discrete groups. Suppose that $A$ is an
 Abelian topological group and $\kappa$ is an infinite cardinal. Then the 
 following are equivalent.
\begin{Enumerate}
\item $A$ is isomorphic to a subgroup of the direct product of $\kappa$ or less
 discrete groups each of size at most $\kappa$.
\item $A$ is isomorphic to a subgroup of
\[\prod_\kappa\left(\bigoplus_\kappa\Q\oplus
 \bigoplus_{p:\rm prime}\bigoplus_\kappa\Z(p^\infty)\right)_d,\]
where $\Z(p^\infty)$ is the quasicyclic group and $\cdot_d$ denotes the
 topological group endowed with the discrete topology.
\item $A$ is TNA and $w(A)\leq\kappa$.
\item $A$ is TNA and $c(A),\chi(A)\leq\kappa$.
\end{Enumerate}
\end{Enumerate}
\end{theorem}
\begin{proof}
(1): We proceed as in the argument for the universality of $\Sym(\omega)$ among
the TNA Polish groups (cf.~\cite[Theorem 1.5.1]{BK}). Since $\Sym(\kappa)$ is
TNA and of weight $\kappa$, we have that $\rm (a)\Rightarrow (b)$. The
implication $\rm (b)\Rightarrow (c)$ is straightforward.

We show that $\rm (c)\Rightarrow (a)$. Suppose that $G$ is a TNA group. Let
$\cal V$ be a local base for $1$ consisting of open subgroups. We may assume
that $|{\cal V}|\leq\chi(G)$. We denote by $\Omega$ the disjoint union of
quotient spaces
\[ \bigcup\{ G/H : H\in \cal V \}.\]
The natural action of $G$ on $\Omega$ such that
\[ g\cdot(xH)=(gx)H \]
induces an isomorphic embedding $G\rightarrow\Sym(\Omega)$. Since
$|G/H|\leq c(G)$ for each open subgroup $H$, we get $|\Omega|\leq c(G)\chi(G)$.

(2): The implication $\rm (a)\Rightarrow (b)$ follows from the universality of
\[\bigoplus_\kappa\Q\oplus\bigoplus_{p:\rm prime}\bigoplus_\kappa\Z(p^\infty)
\]
among (discrete) Abelian groups of size $\kappa$ (see \cite{KS}, Theorem 0.1).

Let $A$ be an Abelian TNA group. By virtue of
\cite[III, \S 7, $\rm N^o$ 3, Proposition 5]{Bo}, there exists an isomorphic
embedding from $A$ into
\[ \varprojlim\{A/B:\mbox{$B$ is an open subgroup of $A$}\}.\]
The estimation of the cardinality is similar as in (1).
\end{proof}

\begin{remark}\label{nam}\rm
A metric $d$ is said to be non-Archimedean if it satisfies the strong
triangle inequality
\[d(x,y)\leq\max\{d(x,z),d(z,y)\}.\]
It is easily seen that a topological group with a right (or left) invariant
non-Archimedean metric is TNA. For metrizable groups, these two notions
coincides: if $\{U_n:n<\omega\}$ is a neighborhood basis of $1$ consisting of
decreasing sequence of open subgroups with $U_0$ the whole group, then
$d(x,y)=2^{-\max\{n<\omega:xy^{-1}\in U_n\}}$ for $x\neq y$ determines a right
invariant compatible non-Archimedean metric.

G.~Rangan asserted in \cite{Ra}: ``Lemma 3.1: Suppose a topological group $G$
is such that its topology is given by a non-archimedean metric $d$ then there
is an equivalent non-archimedean right (or left) invariant metric $\rho$ on
$G$.'' and ``Theorem 3.2: Let $G$ be a separable totally disconnected ordered
topological group. Then $G$ is non-archimedean metrizable.'' (A topological
group is said to be non-archimedean metrizable if there exists a
non-archimedean right (or left) invariant metric on $G$ which induces the
topology of $G$.) 

The additive group of rational numbers $\Q$ is, however, a counterexample to
both of the statements, which is separable, metrizable, zero-dimensional and
TA.
\end{remark}

\section{Disconnected Polish groups}

Among the properties in Remark \ref{dc} (2), no implication other than
indicated holds even for Abelian Polish groups. We construct witnessing
examples from sequence spaces.

\begin{example}\label{ss}\rm
As usual, let $c_0$ be the Banach space of the real sequences converging to $0$
with the norm $||\cdot||_\infty$ and $l^1$ that of absolutely summable ones
with $||\cdot||_1$. We set
\[ R = \{ a\in{}^\omega\R : (\forall n<\omega)
 (a(n)\in(1/(n+1)!)\Z) \}, \]
\[ S = \left\{ a\in\bigoplus_\omega\Z : \sum_{n<\omega}a(n)=0 \right\}, \]
\[ \Gamma_0 = c_0 \cap R, \]
\[\Gamma_1 = l^1 \cap R, \]
where $\bigoplus_\omega\Z$ is identified with
$\{a\in{}^\omega\Z:
 \mbox{$a(n)=0$ for all but finitely many $n$}\}\subset{}^\omega\R$.
Then $\Gamma_0,\Gamma_1,\Gamma_0/S$ and $\Gamma_1/S$ are Polish groups with the
following properties.
\begin{center}
\begin{tabular}{c|cccc}
                     & $\Gamma_0$ & $\Gamma_0/S$ & $\Gamma_1$ & $\Gamma_1/S$ \\
\hline
totally-disconnected & Yes        & Yes          & Yes        & Yes \\
zero-dimensional     & Yes        & Yes          & No         & No   \\
SMOG                 & Yes        & No           & Yes        & No   \\
TNA                  & No         & No           & No         & No   \\
TA                   & No         & No           & No         & No
\end{tabular}
\end{center}
\end{example}
\begin{proof}
$\Gamma_0$ is zero-dimensional: Every open ball is clopen.

$\Gamma_0$ and $\Gamma_1$ have SMOG: Suppose that $a\neq 0$. Then there is $n$
with $a(n)\neq 0$. The open subgroup $\{ b:b(n)=0 \}$ excludes $a$.

$\Gamma_0$ is not TNA: For each $a\neq 0$, we have that
$\lim_{n\rightarrow\infty}||na||_\infty=\infty$. Accordingly every nontrivial
subgroup is unbounded.

$\Gamma_0/S$ is zero-dimensional: Since $S$ is discrete, zero-dimensionality of
$\Gamma_0$ is preserved under the quotient.

$\Gamma_0/S$ and $\Gamma_1/S$ do not have SMOG: Let $e_n$ denote the $n$-th
unit vector $(\underbrace{0,...,0}_{n},1,0,...)$. We show that
$e_0\bmod S(\neq 0)$ belongs to all open subgroups. Let $A$ be an open
subgroup. Then for sufficiently large $n$, we have $e_n/(n+1)\bmod S\in A$.
Accordingly $e_0\bmod S = (n+1)(e_n/(n+1)\bmod S) \in A$.

$\Gamma_0/S$ and $\Gamma_1/S$ are not TA: the subgroup
$\{ b:b(1)\in\Z \}\bmod S$ is open and proper.

$\Gamma_1$ is not zero-dimensional: Theorem \ref{non0} below.

$\Gamma_1/S$ is totally-disconnected: Suppose that $a\bmod S \neq 0$. Then for
$r>0$ sufficiently small, $(B_0(0,r)\cap \Gamma_1)\bmod S$ is clopen and does
not include $a$, where $B_0(a,r)$ is the open ball of $c_0$ with center $a$ and
radius $r$.

$\Gamma_1/S$ is not zero-dimensional: Non-zero-dimensionality is also
hereditary to quotients by a discrete subgroup.
\end{proof}

These examples behave wildly as to other properties.

\begin{example}\label{3smog}\rm
While $\Gamma_0/S$ does not have SMOG, it has a subgroup
$\bigoplus_\omega\Z/S\cong\Z$ with quotient
$(\Gamma_0/S)/(\bigoplus_\omega\Z/S)\cong\Gamma_0/\bigoplus_\omega\Z$ such that
both have SMOG.
\end{example}

\begin{remark}\rm
Since every totally-disconnected locally compact group is TNA, its quotient is
always totally-disconnected as well. On the other hand, as seen from
\[\Gamma_1\left/\left\{a\in \Gamma_1:\sum_{n<\omega}a(n)=0\right\}\right.
 \cong\R,\]
a quotient of a totally-disconnected Polish group may be nontrivially
connected.

Let $c$ denote the space of convergent sequences with the norm
$||\cdot||_\infty$. Then we have $(c\cap R)/(c_0\cap R)\cong\R$. But $c\cap R$,
which has SMOG similarly as $\Gamma_0$, is not zero-dimensional by
Theorem \ref{non0}.
\end{remark}

\begin{prob}\label{0q}\rm
Concerning Polish groups, is zero-dimensionality inherited by quotients?
\end{prob}

We do not know whether the examples in Remark \ref{dc} (4) (5) may be
superseded by a ``strong'' one.

\begin{prob}\rm
Does there exist a zero-dimensional TA Polish group?
\end{prob}

Some ``naive constructions'' toward answering the problems have failed due to
the limitation below for zero-dimensional complete metric groups, which is
obtained by modifying the argument by Erd\H os \cite{Er} for a subgroup of the
Hilbert space $l^2$.

\begin{theorem}\label{non0}
Suppose that $A$ is an Abelian topological group topologized with a complete
invariant metric $d$ and $\Pair{a_n:n<\omega}$ is a sequence in $A$ converging
to $0$. We denote by $C$ and $\tilde C$ the sets of finite subsums and
converging ones, respectively, $\sum_{n\in X}a_n$ for $X\subset\omega$. We set
$\rho=\sup\{d(y,0):y\in A\}\in\R\cup\{\infty\}$ and assume that there exists a
function $\nu:C\rightarrow\R_{\geq 0}$ such that
\begin{equation}\label{div}
\sum_{n=0}^\infty\nu(a_n)=\infty,
\end{equation}
\begin{equation}\label{change}
\mbox{$\displaystyle \sum_{n\in X}\nu(a_n) =\nu\left(\sum_{n\in X}a_n\right)$
 for each finite $X\subset \omega$}.
\end{equation}
\begin{equation}\label{conv}
(\forall\varepsilon>0)(\exists\delta>0)(\forall x\in C)
 (\nu(x)<\delta\Rightarrow d(0,x)<\varepsilon),
\end{equation}
\begin{equation}\label{bound}
(\forall\varepsilon<\rho)(\exists\delta<\infty)
 (\forall x\in C)(\nu(x)>\delta\Rightarrow d(0,x)>\varepsilon),
\end{equation}
Then $\tilde C$ is not zero-dimensional.
\end{theorem}
\begin{proof}
Assuming that $U$ is an open neighborhood of $0$ contained in some open ball
$B(0,r)$ with center $0$ and radius $r<\rho$, we show that $U\cap \tilde C$ is
not closed in $\tilde C$.

We define sequences $\Pair{c_m:m<\omega}$ in $U$ and $\Pair{n_m:m<\omega}$ and
$\Pair{n'_m:1\leq m<\omega}$ of integers by induction such that
$c_m=\sum_{k=1}^m\sum_{n=n'_k}^{n_k}a_n$. We set $c_0=0$ and $n_0=-1$. Since
$U$ is open and $c_m+a_n\rightarrow c_m$ as $n\rightarrow\infty$, there is a
natural number $n'_{m+1}>n_m$ with $c_m+a_{n'_{m+1}}\in U$. Due to (\ref{div}),
$\sum_{n=n'_{m+1}}^\infty\nu(a_n)=\infty$ as well, and hence for sufficiently
large $i$, it occurs that $c_m+\sum_{n=n'_{m+1}}^i a_n\not\in B(0,r)$ by
(\ref{bound}), where we use (\ref{change}) freely. So we may find $n_{m+1}$
such that
\[c_m+\sum_{n=n'_{m+1}}^{n_{m+1}}a_n\in U,\]
\[c_m+\sum_{n=n'_{m+1}}^{n_{m+1}+1}a_n\not\in U.\]

Then we have $c_m\in C\cap U$ and $c_m+a_{n_m+1}\in C\setminus U$ for each $m$.
We show that the sequence $\Pair{c_m:m<\omega}$ is convergent. Since it is in
$B(0,r)$, the sequence of
$\nu(c_m)=\sum_{k=1}^m\nu\left(\sum_{n=n'_k}^{n_k}a_n\right)$ is bounded due to
(\ref{bound}). The latter is increasing as well, and hence it is convergent.
Therefore for any $\delta>0$, sufficiently large $i$ and every $j>i$ satisfy
$\nu(c_j-c_i) =\sum_{k=i+1}^j\nu\left(\sum_{n=n'_k}^{n_k}a_n\right)<\delta$.
Since $\delta$ is arbitrary, (\ref{conv}) yields that the sequence
$\Pair{c_m:m<\omega}$ is Cauchy, so convergent.

Since the sequences $\Pair{c_m:m<\omega}$ in $C\cap U$ and
$\Pair{c_m+a_{n_m+1}:m<\omega}$ in $C\setminus U$ converge to the same point in
$\tilde C$, we conclude that $U\cap\tilde C$ is not clopen in $\tilde C$.
\end{proof}

We exhibit some examples to which the theorem is applicable.

\begin{remark}\rm
The totally-disconnected TA Polish group in \cite{St} is seen not to be
zero-dimensional as follows. If we set $a_n=2^{-2n}$ and $\nu=D$, then the
premises of the theorem hold.

For $\Gamma_1$ in Example \ref{ss}, $a_n=e_n/(n+1)$ and $\nu=||\cdot||_1$ will
do.

As to $c\cap\Gamma_0$, if we set $a_n=(\underbrace{0,...,0}_{n},1,1,...)/(n+1)$
and $\nu=||\cdot||_\infty$, then it goes well. The last argument is based on
the proof provided by Professor Katsuya Eda to the author that the set of
rational points in $c$ is not zero-dimensional.
\end{remark}

\noindent Eda Laboratory\\
School of Science and Engineering\\
Waseda University\\
Shinjuku-ku, Tokyo\\
169-8555, Japan\\
E-mail: higasik@logic.info.waseda.ac.jp

\end{document}